# Seepage Analysis and Control of the Sahand Rockfill Dam Drainage using Instrumental Data

Parvaneh Nikrou[a] , Sajjad Pirboudaghi[b]

[a] *Department of Geography, The University of Alabama, Tuscaloosa, AL 35487, USA*
*pnikrou@crimson.ua.edu*

[b] *Assistant professor, Engineering faculty of Khoy, Urmia university of technology, Urmia, Iran*
*S.pirboudaghi@uut.ac.ir*

## ABSTRACT

The finite element method is an effective numerical technique for the accurate analysis of seepage, enabling the determination of outlet flow and pore water pressures at any point in the dam body and foundation. In this study, the seepage behavior of the Sahand Dam was analyzed using the finite element method and PLAXIS software, and the results were validated with instrumental data. After validation, permeability coefficients were employed to investigate factors affecting seepage at the reservoir's normal water level, such as the cutoff wall, upstream concrete cover, clay blanket, and foundation depth. The results show that the combined implementation of a clay blanket and an upstream concrete cover causes a significant reduction in dam discharge. Finally, by evaluating the effective parameters influencing dam seepage, an optimized model was developed. In this optimized model, the discharge was reduced to one-third of the corresponding value in the initial dam model.

## 1. Introduction

Seismic loads can induce significant lateral forces and deformations in structures, jeopardizing their stability and structural integrity(1). For earthen dams, managing leakage from the body, foundations, and supports is a critical consideration in their design and construction. The issue of leakage from the body, foundations, and supports is one of the most critical criteria in the design and construction of earthen dams(2). As water is stored behind the dam and its level rises, the potential energy of the water particles increases, and it begins to travel in the porous soil. Depending on their size and strength, all dams leak. In fact, all earthen and stone dams are vulnerable to seepage generated by the embankment, foundation, and supports, which wastes water, raises lifting pressure, causes erosion over time, and ultimately jeopardizes the dam's stability(3). This drainage is not a problem; it may not affect the dam's stability. However, difficulties can arise when the rate of seepage unexpectedly accelerates or becomes uncontrollable(4-6). As a result, immediate action is required.

Proper material selection and construction practices are essential not only for minimizing structural deterioration and leakage in hydraulic structures but also for addressing the associated environmental impacts that arise during the design and implementation of seepage control measures(7). While structural integrity remains a priority, raising public awareness about hazard risks, along with enhancing project management approaches(8-10) for maintenance, is essential. This involves implementing appropriate management methodologies for inspections in hazardous zones, with measures to manage traffic at access points(11-13), while also considering the environmental consequences of potential spillway failures.

Effective leakage control strategies prevent water loss behind earthen dams and simultaneously mitigate the buildup of pore pressure within the porous media, which can otherwise lead to reduced effective stress, decreased shear strength, slope instability, and internal erosion. Moreover, controlling leakage forces reduces the risk of soil particle displacement, scouring through the embankment or foundation, and potential destabilization of the entire structure, highlighting the need for an integrated strategy that balances both structural stability and environmental sustainability (14).

In this context, one of the most fundamental criteria for dam design is percolation analysis and leakage calculation using numerical approaches based on complexity. The amount of leakage, the pore water pressure at any point of the dam's body and foundation, the amount of hydraulic gradients in different areas of the dam such as the core, and the water exit locations from the body will all be determined by studying the seepage of an earthen dam(15-20). Many researches have been conducted in recent years on the subject of seepage with regard to the factors influencing the amount of leakage from dams using numerical approaches and various software. comprising the effect of injection(21-24), the effect of the waterproofing curtain(25-30), the calculation of the permeability coefficients of the foundation and support of various dams and drains(31-34), clay blanket(35), the effect of the position and number of drainage and boiling wells(36, 37), and so on. All of these studies are based on numerical modeling, but the results are not completely reliable. because no numerical modeling validation has been performed in this research. Utilizing numerical and hydrodynamic models is a widely accepted method in the literature for simulating and predicting seepage behavior in dam structures, as they allow for detailed analysis of subsurface flow dynamics and potential leakage pathways (38-41). This rigorous, data-driven approach to understanding complex phenomena and optimizing performance through detailed analysis is a hallmark of scientific inquiry across diverse fields (2, 42-45)

In this study, the seepage of the Sahand dam was evaluated using a finite element model, and the results were validated against high-precision instrumental data to develop an optimized seepage control strategy and a robust predictive model for minimizing leakage.

## 2. Case Study

Sahand Dam is one of the largest earthen dams in the country, erected on the Qorangu River, one of the tributaries of the Ghazal Ozen River (twelve kilometres long). The reservoir of the storage dam has a maximum and minimum height of 1600 and 1560 meters above sea level, respectively. The dam is earthen with a middle impervious core at a height of 47 meters from the bed and 59 meters from the foundation, the body's soil volume is 1.3 million cubic meters, the useful volume and total reservoir volume are 148 and 165 million cubic meters, respectively. The crown length and width is 450*10 meters, and the annual flow of the river is 149 million cubic meters(46). Figure 1 depicts a view of the dam's body, lake, and operational structures. Figure 2 also shows the deepest section of the dam.

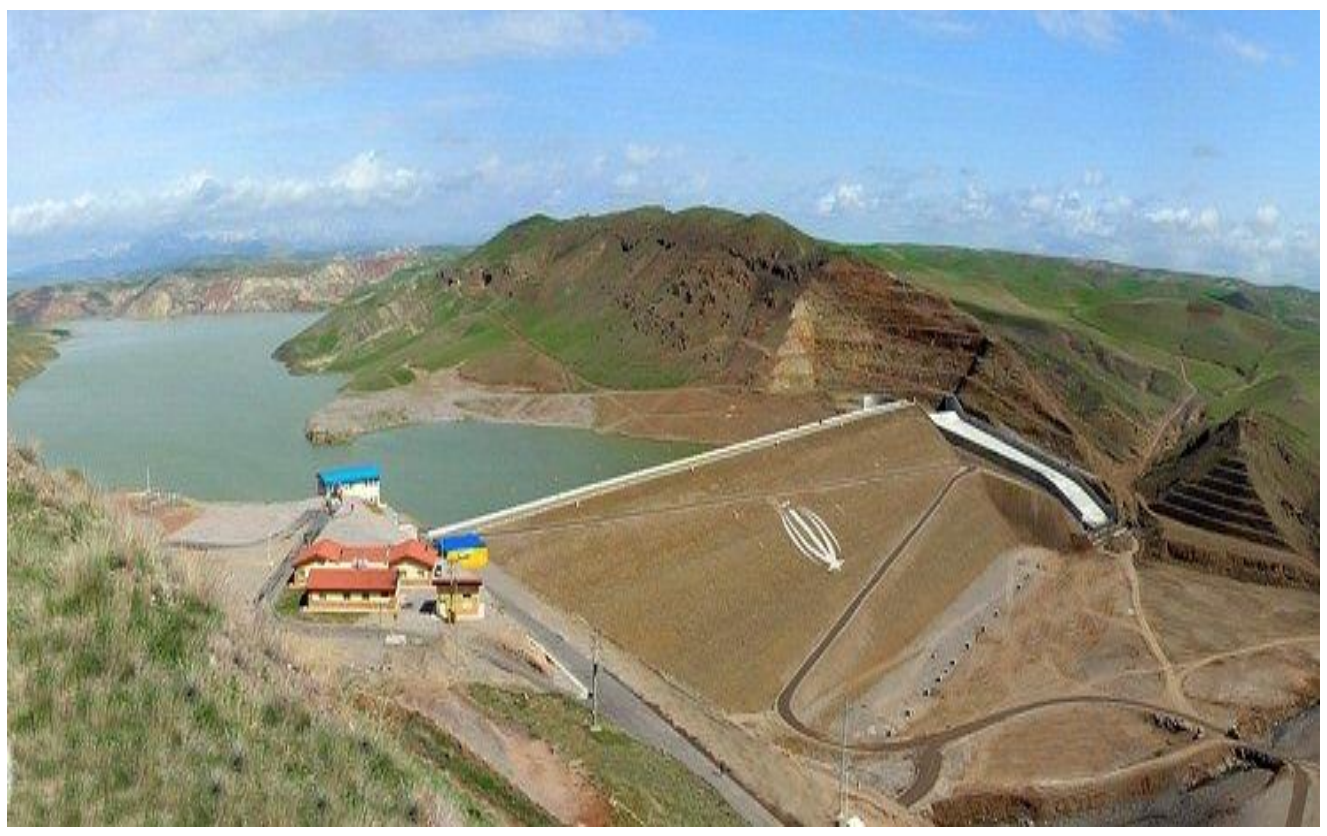

*Figure 1 A view of the dam's body, the lake, and the exploitation buildings(46)*

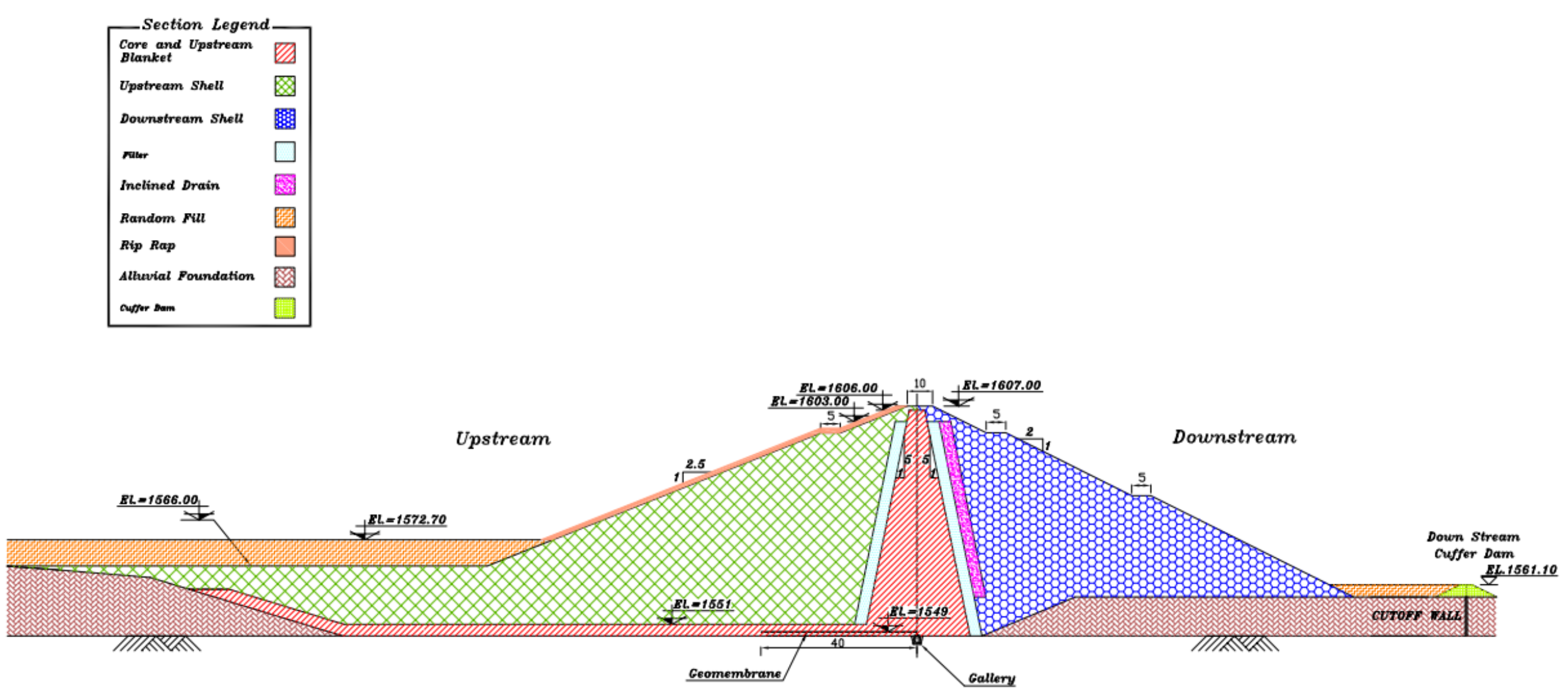

*Figure 2 The deepest section of the dam*

## 3. Seepage Analysis

PLAXIS 8.5 was used to perform seepage studies on the Sahand dam (47, 48). In this analysis, the deepest section of the dam was chosen, and the water level was taken at the maximum level recorded throughout the operation period, i.e. 8/1582 dated 2007/05/07. It should be noted that the tank's regular level is 1600.3. Figure 3 depicts the desired area of the analysis. According to Figure 3, the modeling stages are such that all of the dam's elements, including the asymmetric shape of the stone bed, clay core, filters, and access roads, are represented in geometrical detail. The minimum foundation depth considered for seepage analysis is 60 meters, which is about equal to the dam's hydraulic height. The width of the dam's cross-section is around 300 meters, but the length of the model is 500 meters, therefore it has no effect on the results. Table 1 shows the parameters evaluated for percolation analysis.

The results of the permeability test during drilling and preliminary studies prior to dam building were used for the foundation layers. The values for the core and shell were derived from preliminary loan source studies conducted before to the dam's construction. Because the water must be drained from this component without any pressure, the permeability coefficient of the filter is likewise a significantly bigger number as compared to the other selected values. Due to the dissolution and pressure of the foundation layers, the permeability coefficient of the stone foundation is somewhat greater than before dewatering to provide for more realistic conditions, because the water passage has increased through time. Following the excavation of the dam site, there has been a leakage issue, the majority of which is most likely due to the foundation, raising concerns about the dam's stability. As a result of the leakage flow from the dam's downstream, the water level of the reservoir has not been filled to its normal level in recent years. The adjusted permeability coefficient of the rock foundation is considered to be $10^{-4}$.

Following all of the aforementioned cases, the proposed model has been such that by correcting and correcting multiple models with varying permeability coefficients, it is possible to reach the real conditions that are adapted to the exact dam instruments, i.e. body piezometers. Furthermore, the distribution of the water flow network and pressure graphs should be appropriate and not restricted to a few layers adjacent to each other. The final results are shown in Table 1.

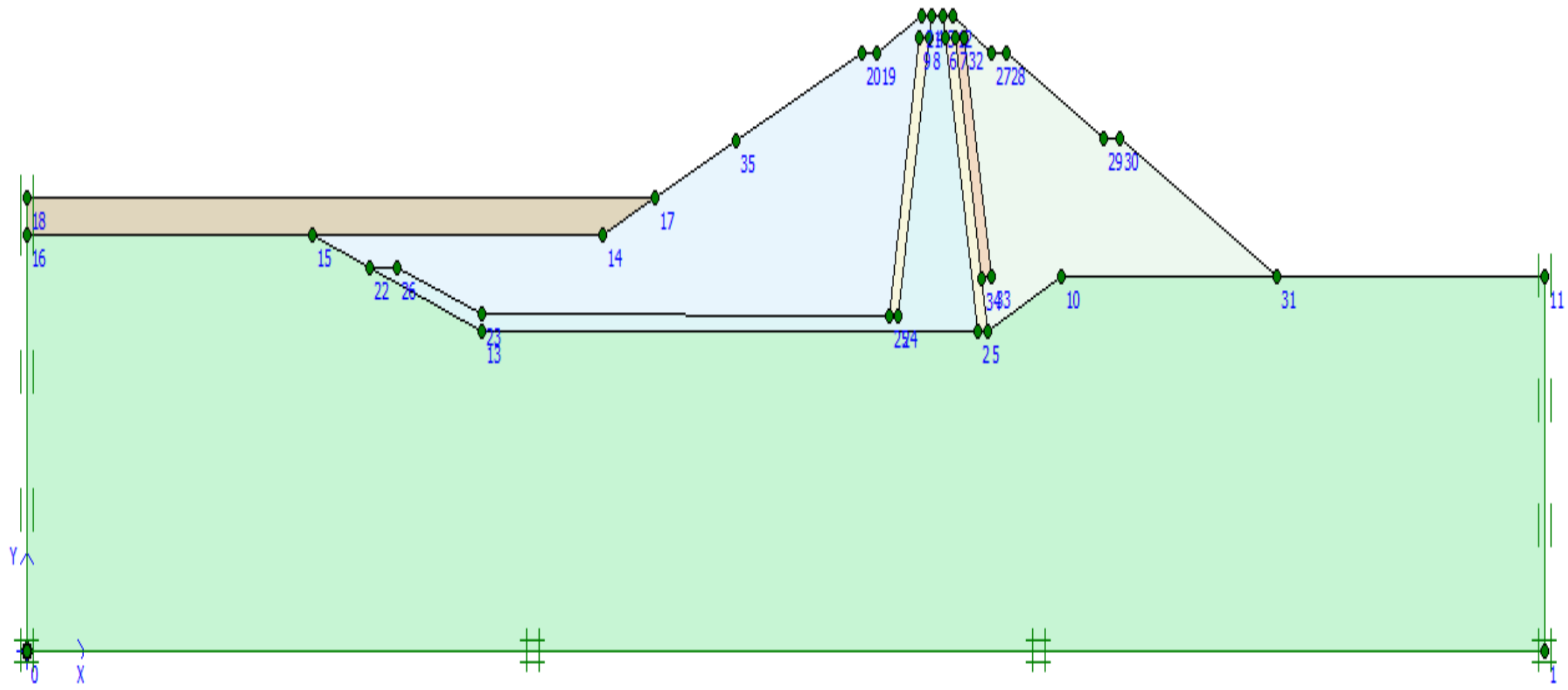

*Figure 3 The section considered for the analysis of seepage of Sahand Dam*

*Table 1 Parameters considered for percolation analysis*

| Material | γ (kN/m³) | φ (°) | C (kN/m²) | K (cm/s) |
|---|---|---|---|---|
| Upstream shell | 20 | 35 | 30 | $10^{-1}$ |
| Downstream shell | 20 | 35 | 30 | $10^{-1}$ |
| Core | 20 | 30 | 50 | $10^{-8}$ |
| Stone Foundation | 21 | 35 | 0 | $10^{-4}$ |
| Filter | 18 | 35 | 0 | $10^{-2}$ |
| Drain adjacent to the filter | 18 | 35 | 0 | 1 |
| Bottom waste material | 19 | 30 | 0 | $10^{-2}$ |

### 3-1. Reservoir water level and instruments

To measure deformation and pore water pressure in the dam's body and bed, accurate instruments are installed and read on a routine basis. The changes in instruments are critical, and with the data gathered, it is able to forecast and control the sealing system's efficiency, the dam's proper functioning, and the failure factors. Figure 4 depicts piezometers of the dam type (middle and deep parts of the dam) installed beneath the body.

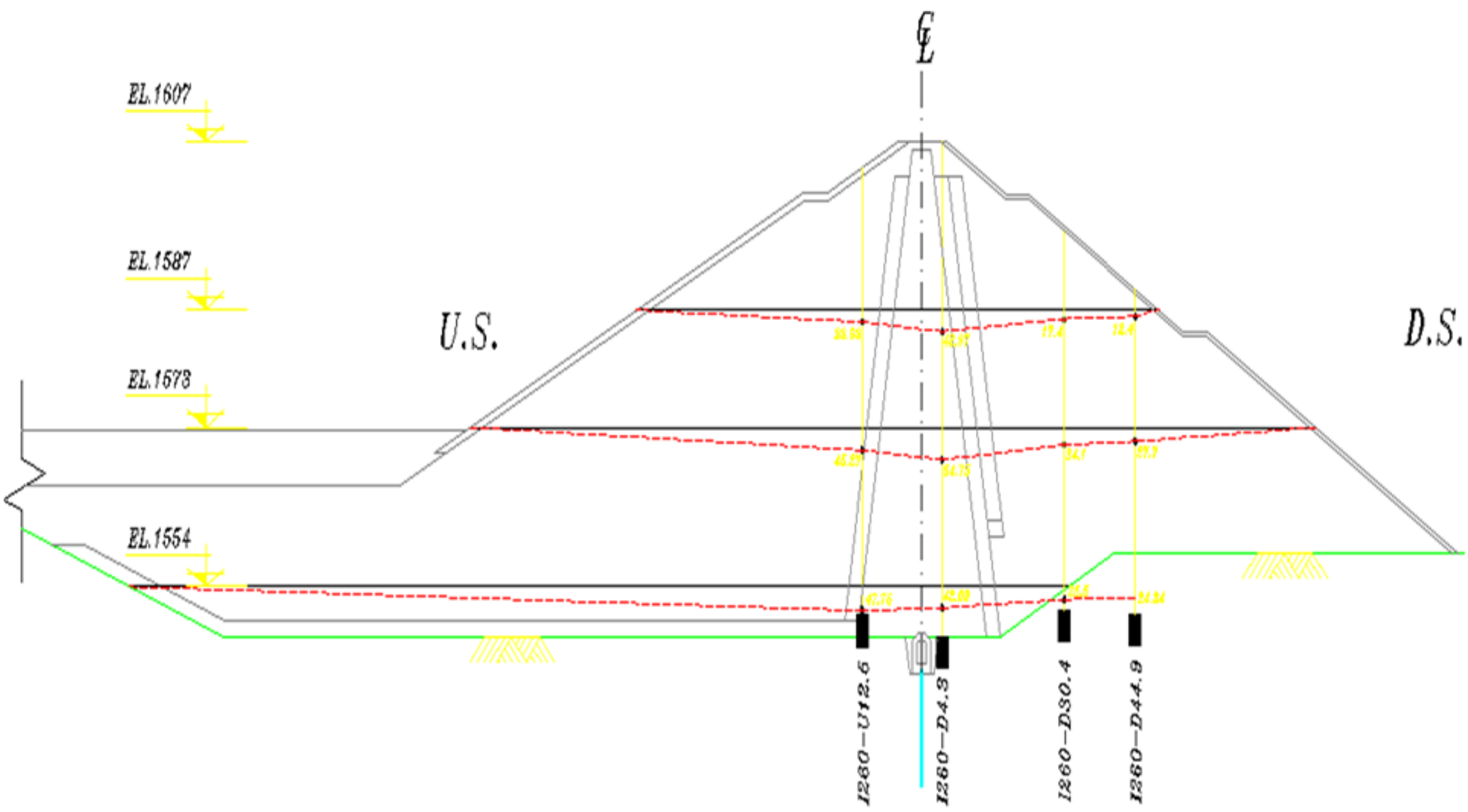

*Figure 4 Cross-sectional piezometers of the dam type*

Figure 5 depicts the variations in the reservoir's water level over time, as well as the changes in the water level of the dam body's piezometers, from the first spring of 2005 to the end of summer 2006. According to Figure 4, these piezometers are installed in various places of the dam. The measured results reveal that the fluctuations in pore water pressure are reasonable throughout the year, and its changes in healthy piezometers are primarily based on changes in tank water level, indicating that the water system is running properly. Of course, in defective piezometers, no readings were made or changes in the water level inside the piezometer did not react to changes in the water level in the tank, indicating piezometer failure and blocking, which requires further behavioural measurement and management.

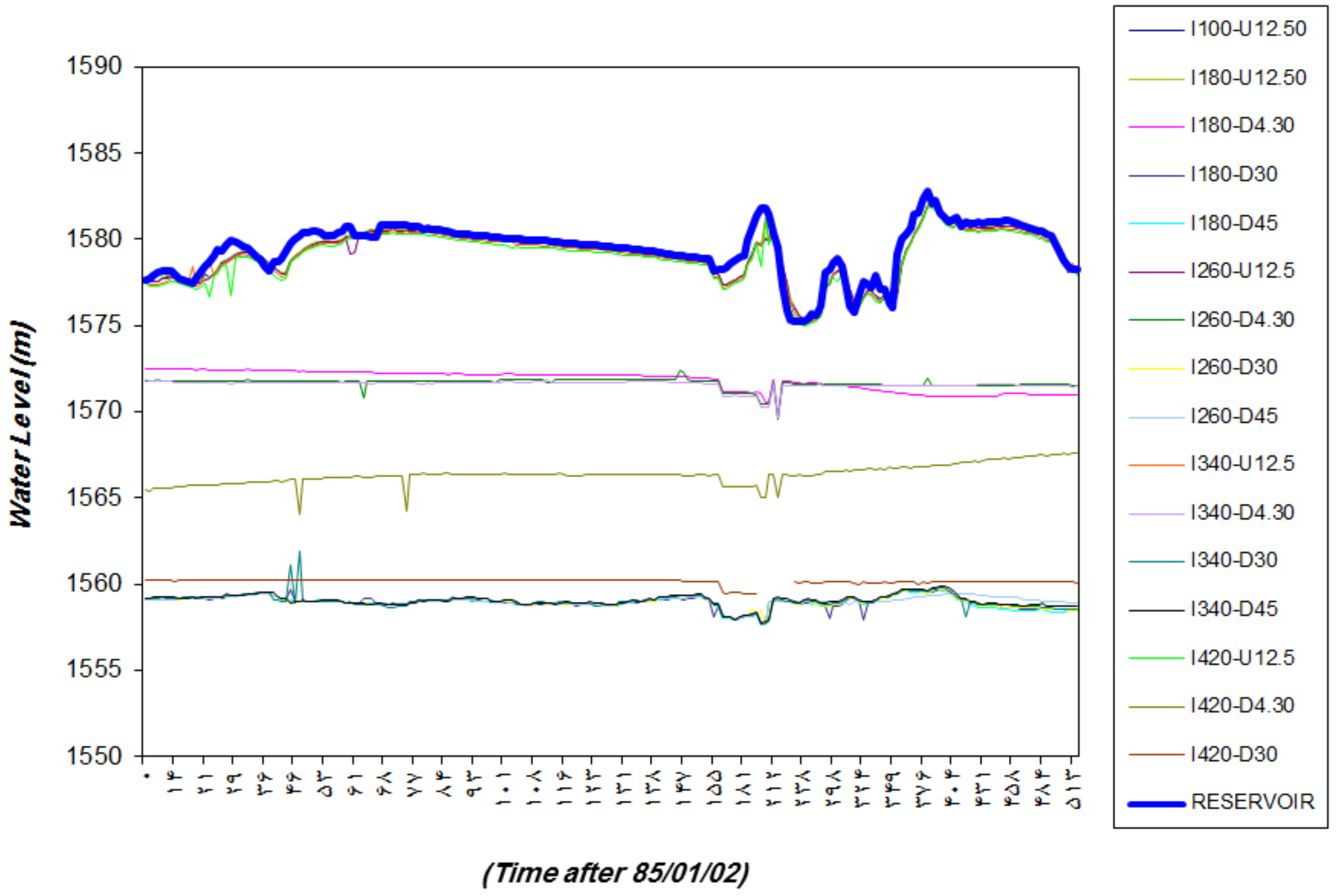

*Figure 5 The readings of piezometers of the dam body compared to the reservoir level*

Most of the piezometers are flawed in the body and do not respond to changes in the tank's water level. Meanwhile, because the most essential section of the dam was chosen for the seepage study, the piezometers of the middle section of the dam depicted in Figure 4 have been used as a basis for validating and anticipating the results of the seepage analysis. Figure 6 illustrates the changes in the flow of leaking water from the downstream valley compared to the reservoir level variations.

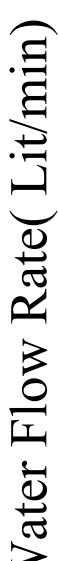

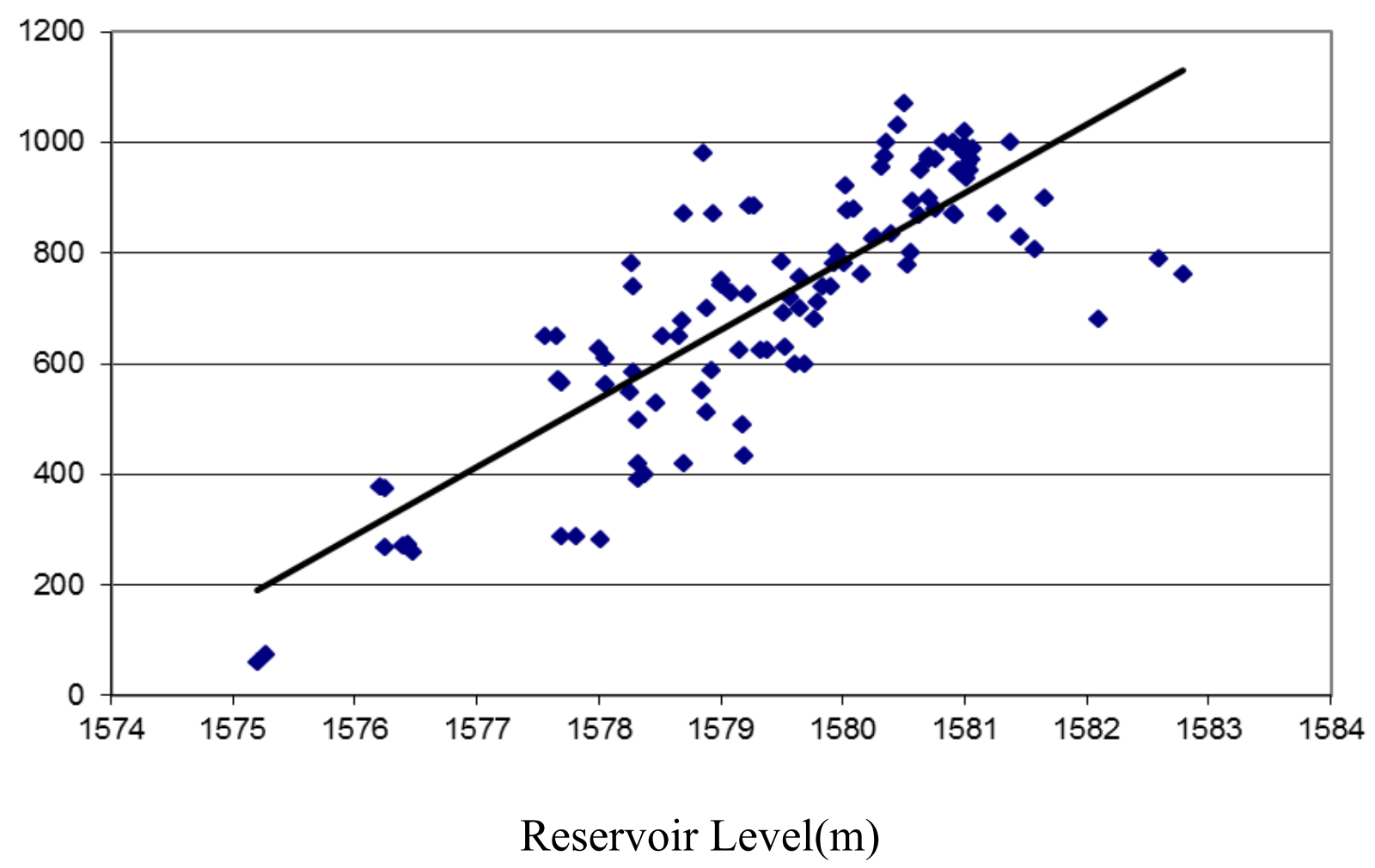

*Figure 6 Changes in the flow rate of the downstream valley's releasing water in relation to reservoir level changes*

### 3-2. Verification Process

Because of their better and more logical position, body-level piezometers in the middle part are used as a basis for validating seepage analysis results. According to the readings of the existing tank's water level, the highest level is associated with the date of 2007/05/07 with the value of 1582/8, which is why it was chosen for modeling with the maximum water load and the results of the readings of four piezometers of the section level. Using the model coordinates in the software and referring to Figure 5 and the available data, water levels for piezometers I260-U12.5, I260-D4.3, I260-D30.4, and I260-44.9 are 93/93, 83/95, 71/47, and 71/3, respectively. It should be noted that the number 260 represents the distance from the beginning of the crown, D and U reflect downstream or upstream, and the following number is the distance of the installation site from the dam section's middle axis. The best comparison between the simulation results of the water flow level inside the body and the water level readings of four body piezometers in the dam's middle section, namely I260-U12.5, I260-D4.3, I260-D30.4, and I260-44.9, is shown in Figure 7. The results are in good agreement even in the area between the filter and the core due to their proximity, indicating that the body material properties, particularly the permeability coefficients and modeling steps, were correct. This kind of dual-input validation—aligning modeled and measured data—echoes principles also supported in other disciplines (see (49)

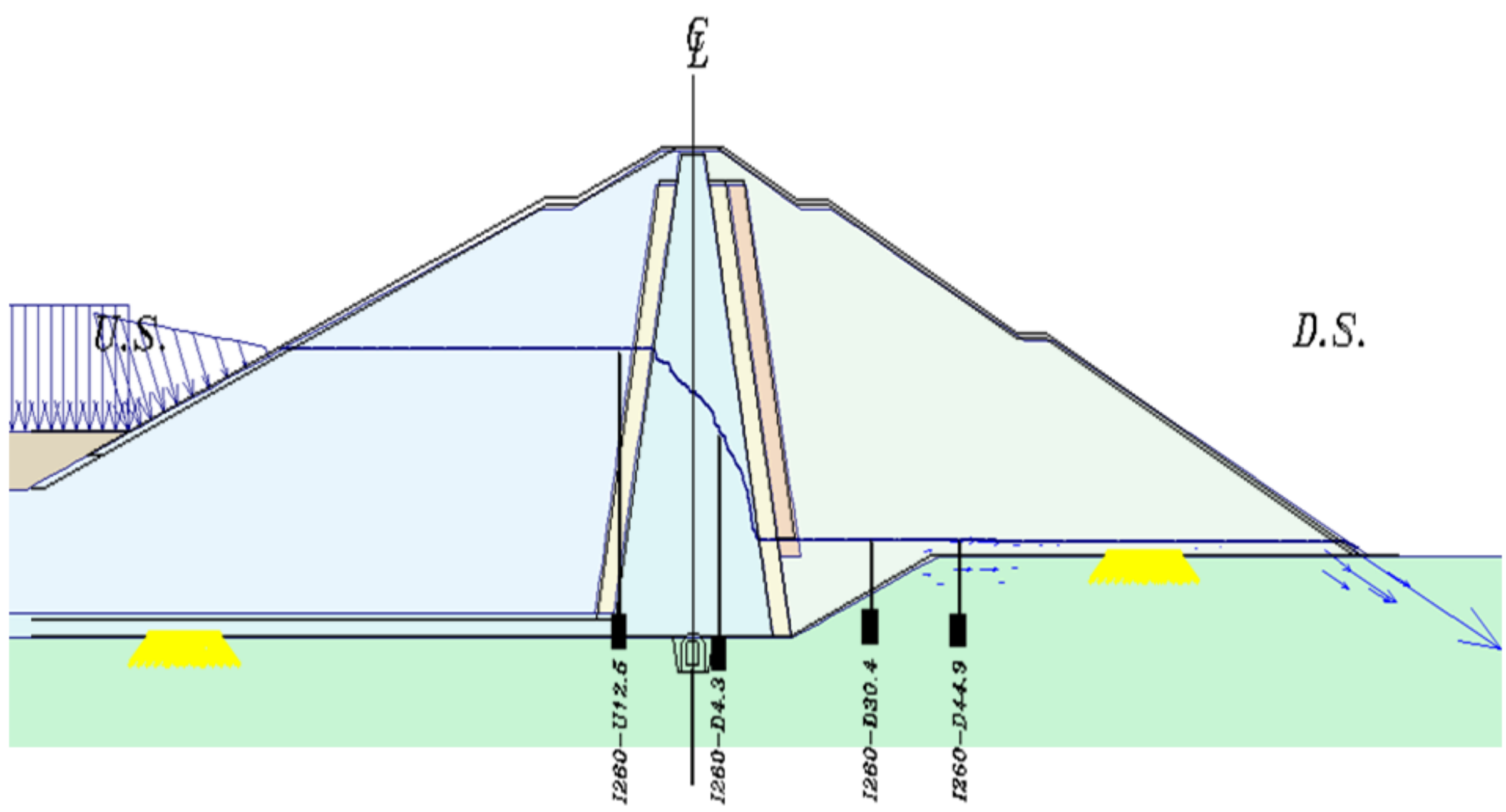

*Figure 7 Comparing the water level of the model and the results of piezometers of the dam body in the middle section*

According to the analytical results, the highest water flow speed inside the body at the level of 1582.8 is 54-10 x 54.5 m/s, and the total seepage volume is 6-10 x 5.5 m3/s per one meter perpendicular to the dam, which Taking into account the length of the dam crown (450 meters), the amount of leaked water is approximately 2.5 litres per second, and if the reservoir water level rises, it will increase up to the normal level. Figure 6 shows that on 2007/05/07, the amount of water leaking from the downstream valley was 760 litres per minute or 12.7 litres per second, which is much higher than the model value of 2.5 litres per second. While body piezometer readings have been used to confirm the accuracy of the created model. As a result, the dam clearly has a leaking problem, and the water level in recent years has not been filled up to the reservoir's normal level, owing to the significant leakage flow from the downstream side of the dam. As a result, by giving solutions, the amount of drainage should be minimized. Also, as shown in Figure 7, the velocity vectors in the path under the downstream shell or the boundary between the rock foundation and the shell and the end part of the dam's outlet are large, which will naturally cause the destruction of these parts, which must be done by reducing hydraulic gradient measures, so that the water can be drained and removed more easily. Because water may wash the shell and develop pores, causing the shell to fail and settle downstream. The simplest remedy appears to be the installation of a ~~table~~ blanket ~~and~~ or claw drain, which can be built of filter material to prevent its effects. As a result, for the analysis of the dam model, a ~~table~~ blanket drain with a depth of 3 meters and a claw with a depth of 10 meters, made of the material of the drain in front of the filter, with a permeability of one centimetre per second, have been used.

### 3-3. Analysis of seepage at normal level

As previously stated, the reservoir water level has not been filled to the normal level by observing and measuring the outflow from the downstream valley in order to prevent the creation of critical conditions owing to high leakage and to ensure the proper operation of the dam. while the dam is designed and built to work at a normal level. As a result, in order to use the dam's full capacity, the analysis of seepage at the tank's normal level, 1600/3, is considered. The results of the previous level, 1582/8, were validated by matching the seepage level of the PLAXIS model with the measurements of the water level of the piezometers in the middle part on the same day and level. It was also demonstrated that the values in **Error! Reference source not found.** were accurately chosen and that the modeling steps were followed correctly. As a result, the previously mentioned parameters and procedures are also used in normal-level modeling. Figure 8 depicts the flow network inside the body at the normal level. The maximum speed of water flow inside the body is 248 x 10-9 m/s, and the total seepage amount is 9.7 x 10-6 m3/s per one meter perpendicular to the dam, which is by considering the length of the dam crest (450 m), the amount of water leaked is approximately 4.4 litres per second, which seems to be much higher than the observed output discharge if the dam is actually filled to the normal level. Therefore, the amount of leakage and output flow should be decreased by offering solutions to achieve normal dam conditions. It can be seen that the continuous deep blanket from the dam's core has performed excellently, forcing the water to enter a reversible state after penetrating the shell.

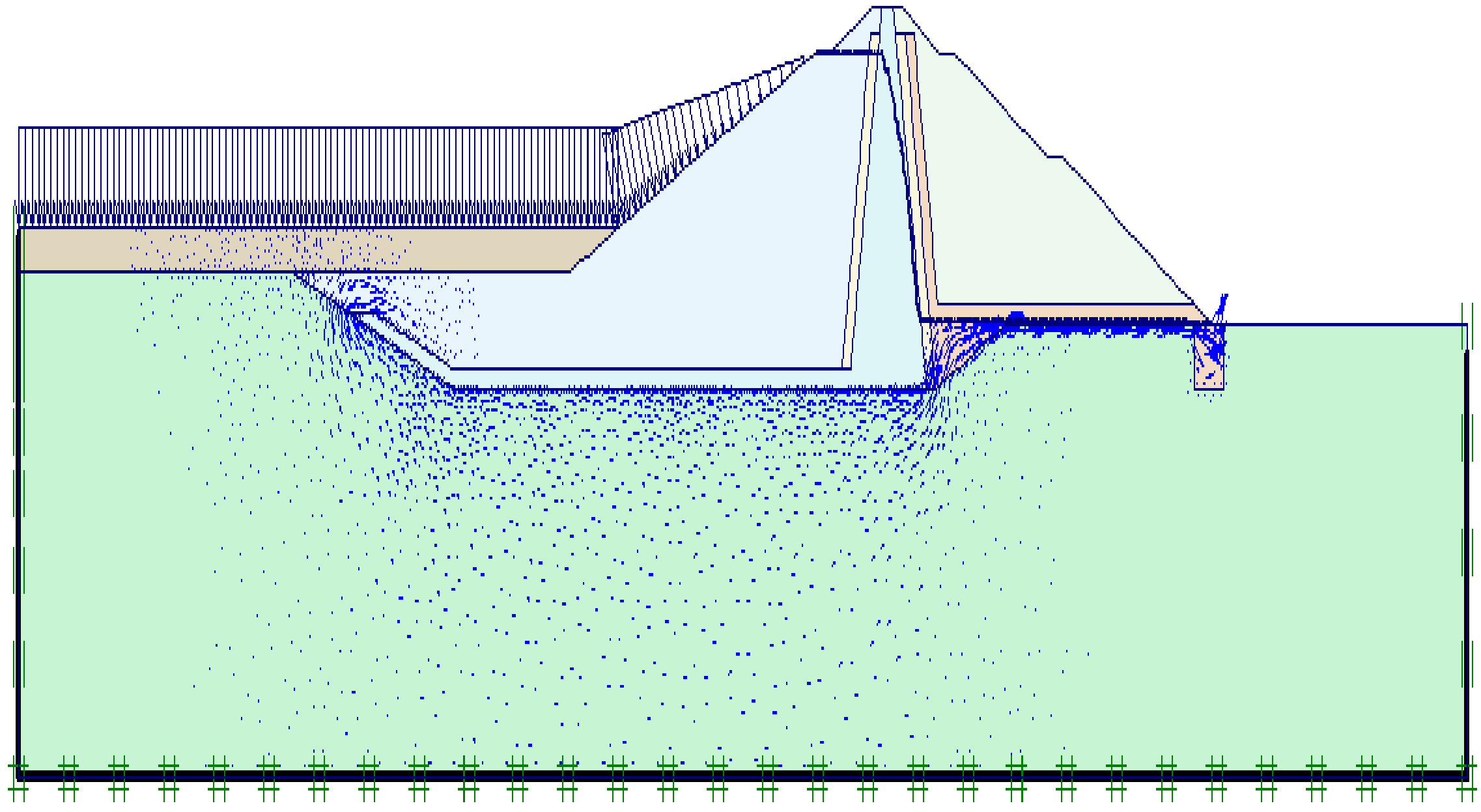

*Figure 8 The flow network at a normal level*

## 4. Investigating solutions to reduce leakage

After predicting and verifying the instrument data, the factors influencing the amount of water leakage from the dam, such as the creation of a watertight curtain, upstream concrete cover instead of rip-rap, the implementation of clay blankets, and the impact of foundation depth, will be investigated separately in the following, so that a composite and optimal model can be obtained from these researches.

### 4-1. The effect of the sealing curtain

After seeing significant dam leakage, the seal wall can be considered under the core. However, after the dam is built, it is not viable to build it below the core, and it must be built in the upstream heel to a depth of, for example, 30 meters. Given the dam's importance and height (59 meters), the clay core under the dam can be sloped to a lower depth and, in some areas, paired with a concrete heel for greater seismic performance. This wall is often constructed of plastic concrete. Figure 9 depicts the results of a seepage analysis using a 2-meter-thick dam wall with the same permeability as the dam core under the core, and Figure 10 displays the results using the wall at the upstream heel. According to the findings, the amount of seepage with the wall under the core is 3.8 litres per second, and the amount of seepage with the wall in the upstream heel is 3.9 litres per second. In fact, the existence of the wall in the upstream heel has a less substantial effect on seepage than the presence of the wall below the core, even though the upstream wall is positioned at a depth of 20 meters inside the shell.

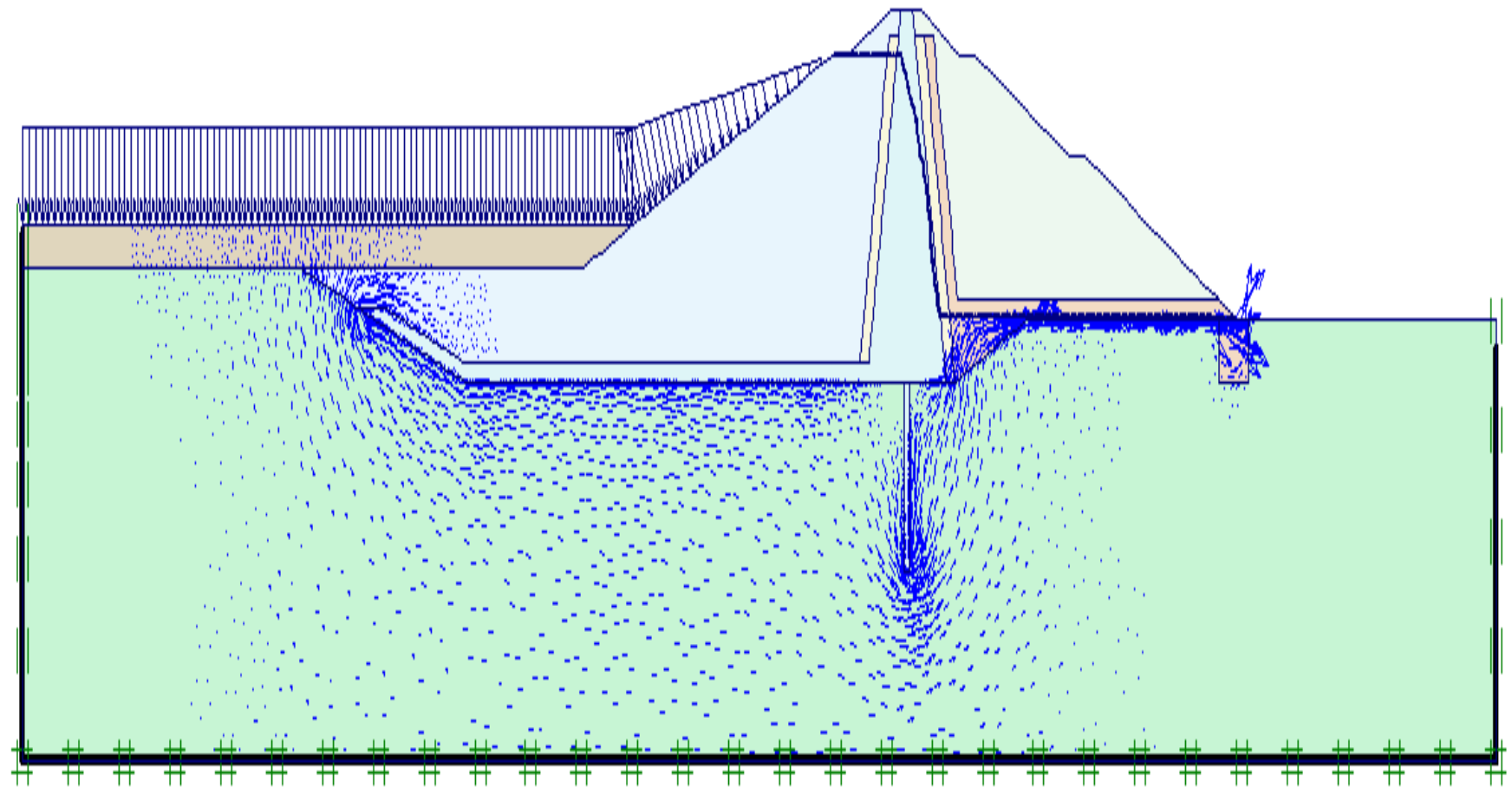

*Figure 9  Flow network with water curtain under the core of the dam*

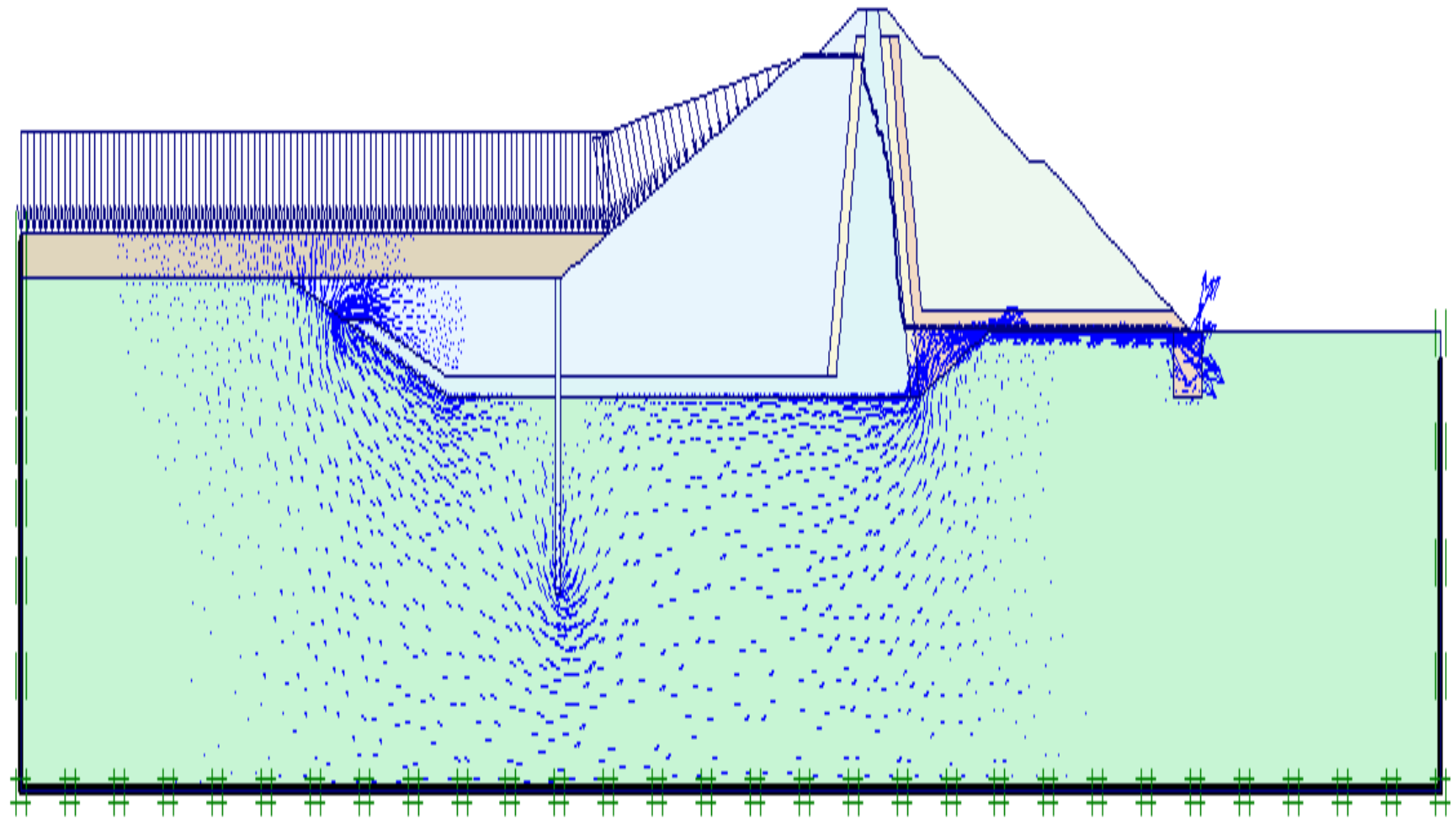

*Figure 10  Flow network with a water curtain at the upstream heel of the dam*

### 4-2. The effect of upstream concrete cover

The upstream slope is covered with riprap. Then, in the following years, when a lot of leakage is observed, it can be turned into an impermeable concrete cover, so that part of the losses caused by the large amount of water seepage can be removed.

The execution process is as follows: initially, this cover is divided into distinct blocks that are molded and concreted, and then the gap between these blocks is sealed with mastic (a type of polymer glue mixed with cement). Figure 11 depicts the results of a seepage analysis with an upstream concrete cover of 0.5 m thickness and permeability of 5-10 cm/s. The amount of leaking water is 4.4 litres per second, based on the seepage results. Contrary to expectations, the percolation from the body and foundation does not change since the water percolates through the foundation layers rather than the shell from the tank's bottom. Figure 11 shows that in the absence of an upstream concrete cover, water tends to enter through the foundation levels. In fact, there is no advantage to building a concrete cover upstream.

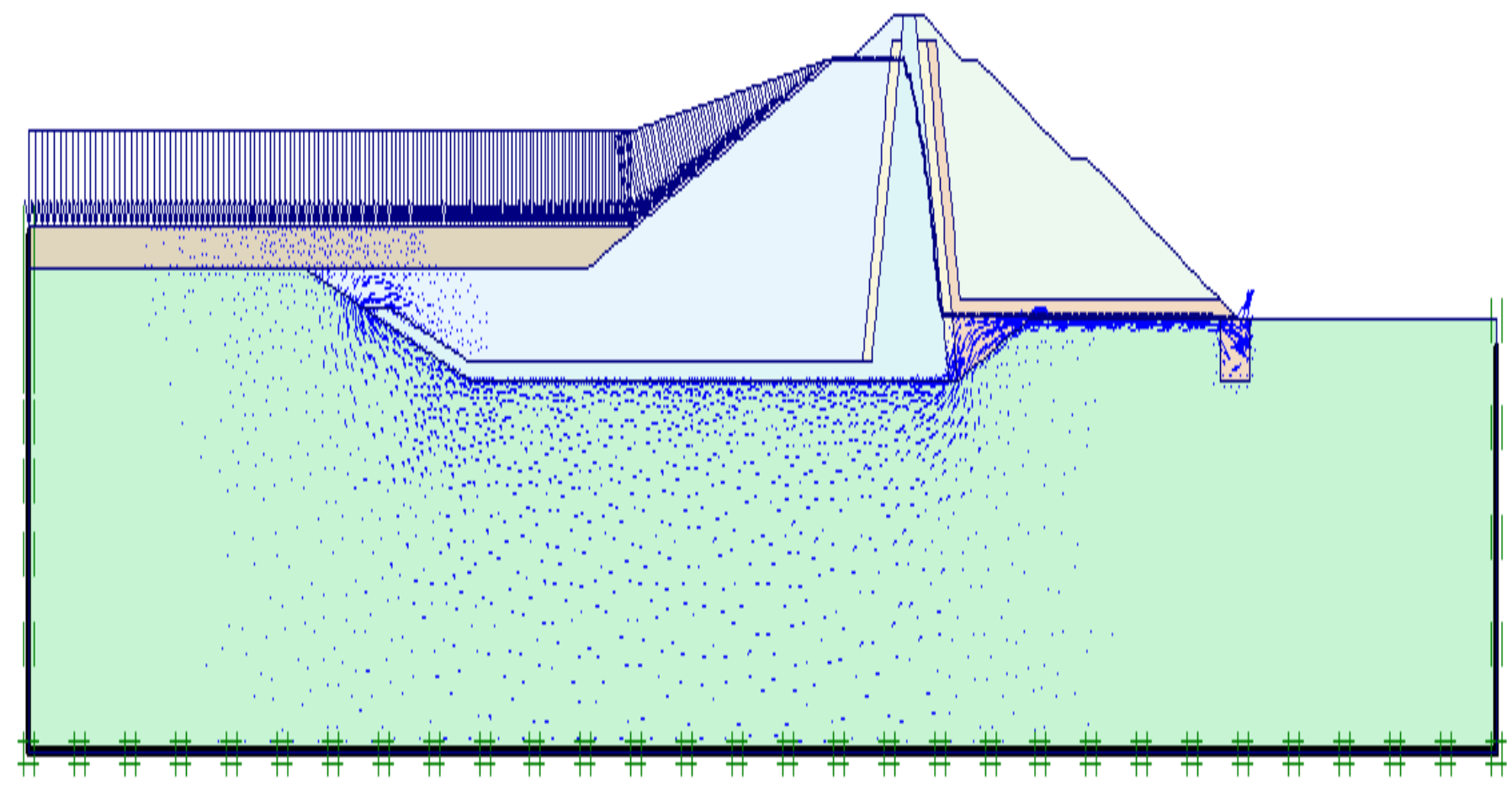

*Figure 11 Flow network under the influence of concrete cover upstream of the dam*

### 4-3. Clay blanket effect

After noticing a lot of leakage and draining the entire water in the tank, a clay blanket can be laid in the bottom of the tank in the shape of an impenetrable layer. However, the downside of the clay blanket is that the over-emptying of the tank is less than the amount required by the farmers, and the drying of the pounded clay causes cracking and loss of the clay cover to some extent (clay blanket). Although this excessive discharge endangers the environment and the animals in the vicinity, it must be resolved by properly coordinating the input flow rate of the tank and the output flow rate of the consuming area. Figure 12 depicts the results of a percolation analysis performed using a 1 meter thick clay blanket at the reservoir's bottom with the same permeability as the 200-meter-long core.

Clay blanket can be applied on the bottom of the reservoir and after draining the water in the reservoir after seeing a lot of leakage in the form of an impermeable layer. However the disadvantage of the clay blanket is that over-emptying the reservoir due to the farmers' need for a level lower than the dead volume and the drying of the pounded clay will to some extent cause cracking and loss of the clay cover (clay blanket).

Using the seepage results, the amount of leaking water is 4.4 litres per second. Contrary to expectations, the percolation from the body and foundation does not alter since the water enters the body from the shell rather than the bottom of the tank and will re-percolate through the foundation layers in a reversible state. As a result, it is advised that the clay blanket be installed simultaneously with the upstream concrete cover, or that a continuous deep blanket of the dam core be installed upstream, completely to the reservoir's bottom.

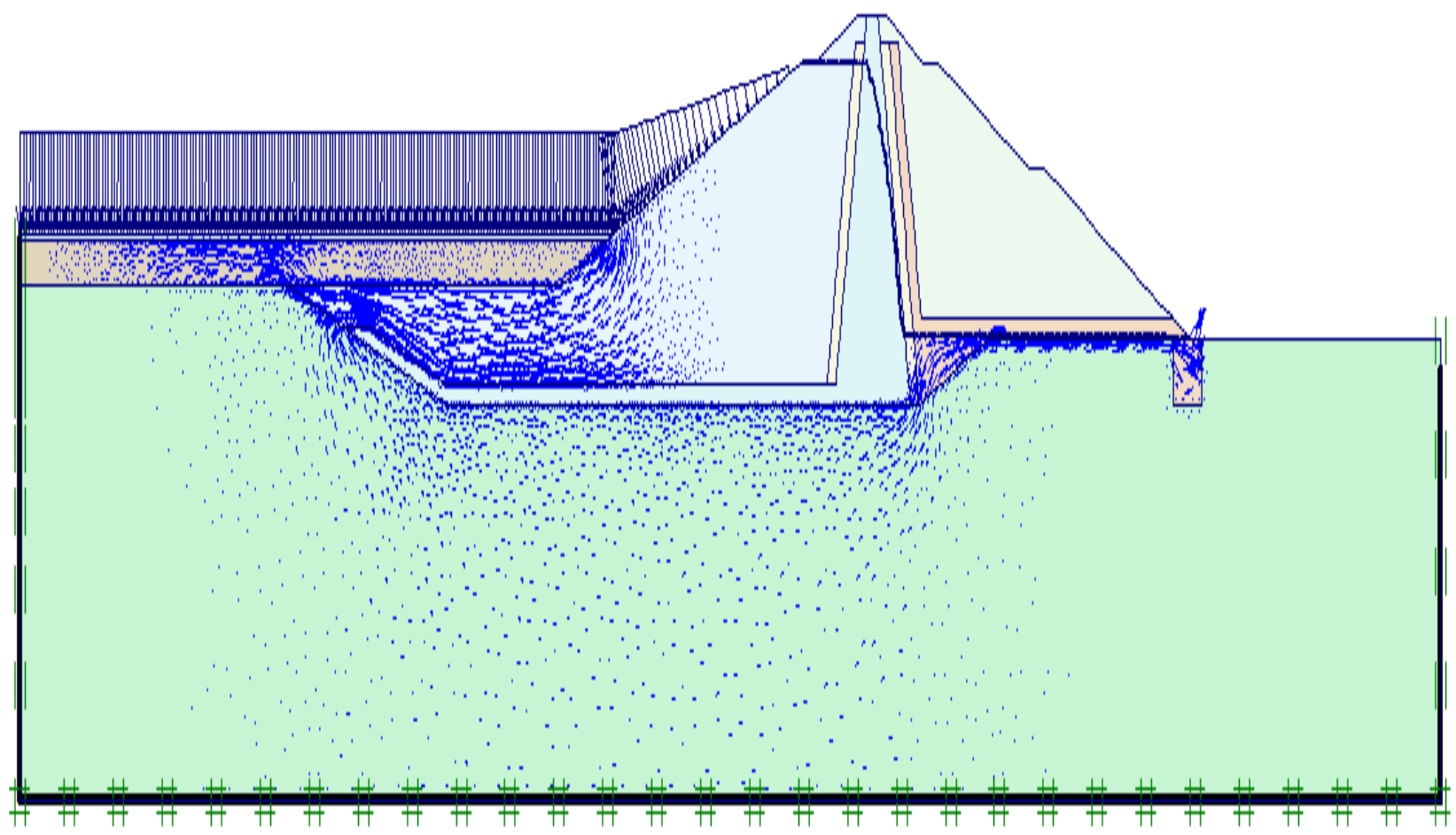

*Figure 12 Flow network with clay blanket on the bottom of the tank*

### 4-4. The effect of foundation depth

In most cases, the depth of the foundation is assumed to be equal to the hydraulic height of the dam during modeling. To explore the influence of foundation depth on dam seepage, three models with depths of 30, 90, and 120 meters were created based on the dam's height (59 meters). Figure 13-15 illustrate the findings of seepage analysis at various depths.

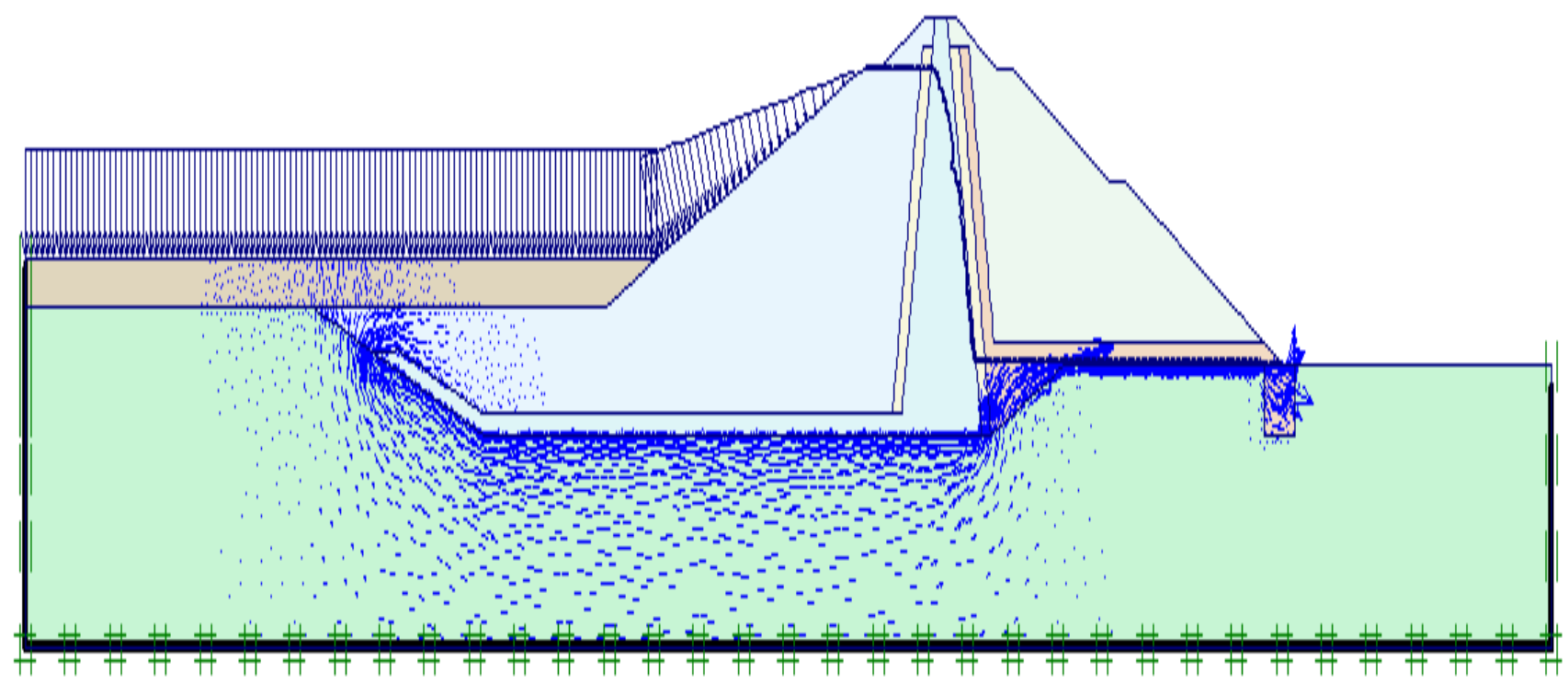

*Figure 13 Stream network with a foundation depth of 30 meters*

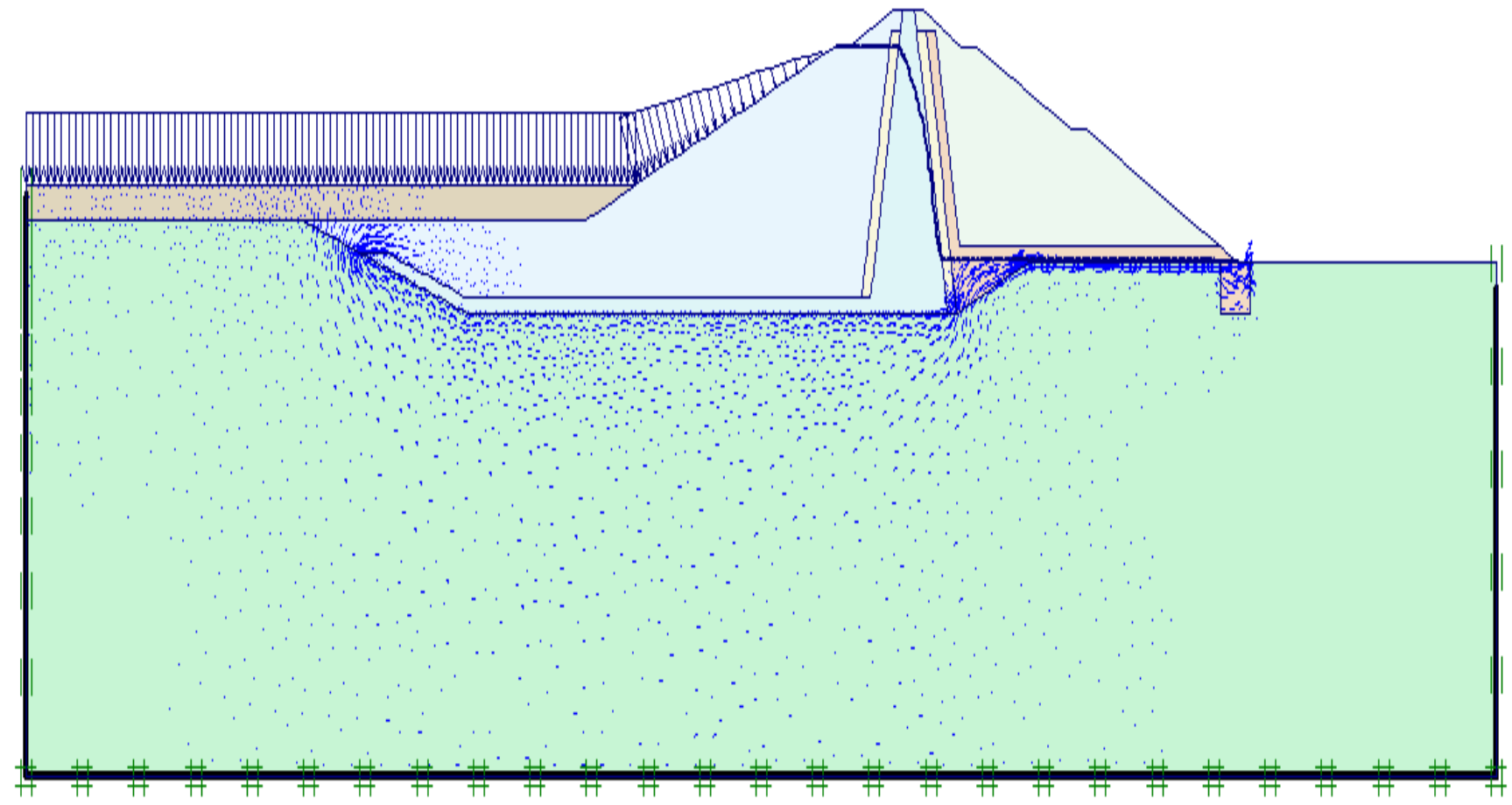

*Figure 14 Stream network with a foundation depth of 90 meters*

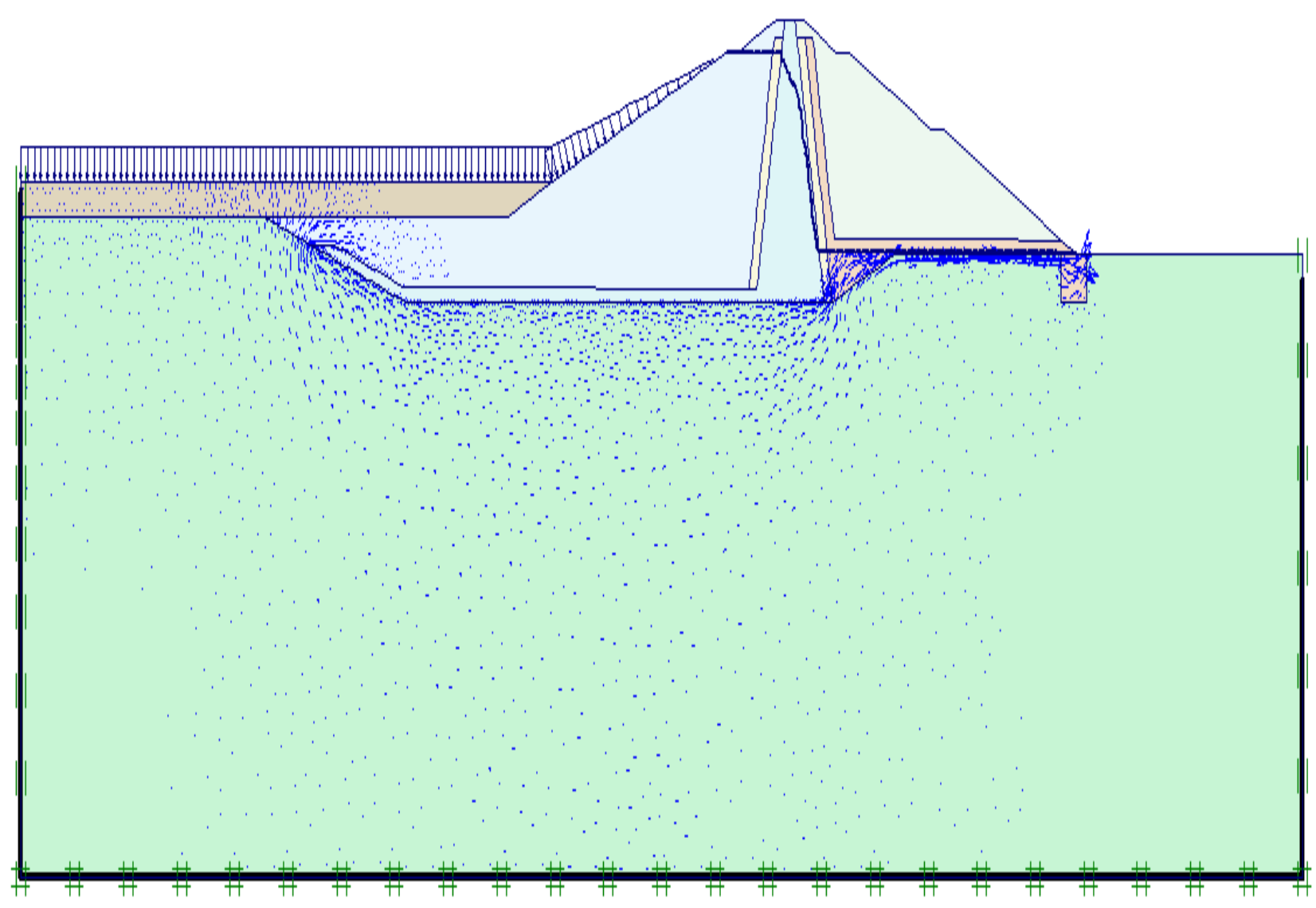

*Figure 15 Stream network with a foundation depth of 120 meters*

Using the seepage data, the amount of leaking water at different depths of 30, 90, and 120 meters is 2.5, 5.9, and 6.8 litres per second, respectively. Furthermore, the seepage value for the normal depth was already 4.4 litres/second. It is obvious that the effect of foundation layers on seepage reduces with increasing depth. According to the results, the foundation modeling is not as suitable as the dam's hydraulic height, and at least 1.5 times the dam's height, i.e. 90 meters, should be used in order to include more realistic dam circumstances for seepage analysis.

### 4-5. Composite optimal model

All the effective components have been incorporated into the dam's final model after evaluating the factors impacting seepage in order to give a full model emerging from this research. The foundation depth is assumed to be 90 meters for this reason. The core's continuous deep blanket is stretched upstream from the incomplete condition of the tank's bottom. The dam wall was built in the upstream heel, and the concrete top cover and clay blanket of the tank bottom were added, along with the previously indicated parameters, before the dam model was examined. Figure 16 depicts the findings of the percolation study and flow network within the body, as well as the foundation for the final model. It should be noted that a finer mesh was used in some areas if necessary.

As can be seen, the likelihood of water penetration from the tank's bottom is limited in the final form. Using the seepage results, the amount of leakage in the final model is less than 1.5 litres per second, which is about one-third of the initial dam model. This is even though the use of the clay blanket in conjunction with the upstream concrete cover had the largest impact on leakage flow and should have been addressed before the dam was built.

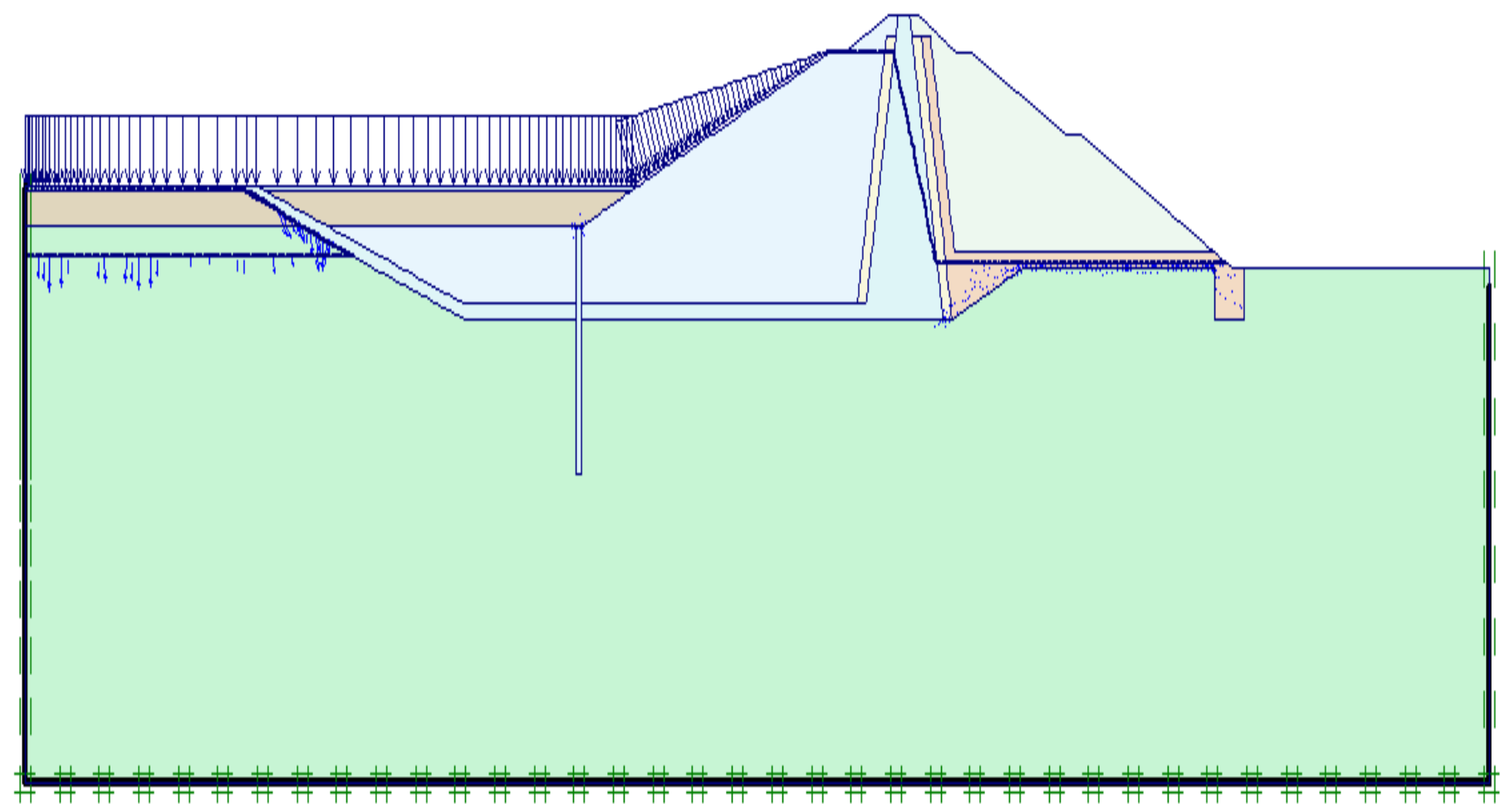

*Figure 16 Flow network and maximum water velocity in the final model*

Given the high gradient around the core and the passage of more infiltrating water from this region, as well as the current lack of a waterproof curtain, it appears that repairing the defective piezometers is important for a more accurate measurement of the dam's behaviour. Because of their position, the piezometers in this section are critical for assessing behaviour, and variations in pore water pressure should be managed in them. An adequate sealing system capable of preventing water infiltration should be installed along the dam's axis to the depth where leakage is possible. It is possible to implement the watertight curtain system in the upstream heel in a deep way along the full length of the dam for this purpose. In addition, the continuous deep blanket of the core should be extended to the bottom of the tank in such a way that the possibility of infiltrated water returning to the upstream side is eliminated, in conjunction with the simultaneous implementation of the clay blanket on the bottom of the tank and the impermeable concrete cover of the upstream surface.

As a potential solution, local injection can be used to lessen percolation by carefully analyzing the geology of the layers, especially dissolvable layers, and identifying their position, if possible. According to statistical research, a green-yellow-red zone should be constructed for the tank before adopting corrective actions to limit the final leak to avoid overfilling the tank and causing critical conditions in the current scenario.

## 5. Conclusion

The current study performed a comprehensive seepage analysis of the Sahand Dam using the Finite Element Method (FEM) and validated the results with instrumental data. After verification, permeability coefficients were employed to evaluate the influence of various seepage control mechanisms, including the reservoir's cutoff wall, upstream concrete cover, clay blanket, and foundation depth on seepage behavior under normal water levels. The key findings from this study are as follows:

- Impact of Table and Claw Drains: The installation of table and claw drains effectively lowers the hydraulic head within the downstream shell, facilitating more efficient water discharge and mitigating the risk of claw damage and scouring.
- Leakage Identification: The comparison between the model's computed outflow and observed downstream valley discharge reveals a leakage issue, as field observations indicate that the reservoir's water level has not reached its designed capacity in recent years. Therefore, appropriate remedial measures should be implemented to reduce the drainage and restore reservoir efficiency.

- Quantification of Leakage: At the reservoir's normal level, the model estimates a seepage rate of approximately 4.4 liters per second. However, the actual discharge is expected to be higher if the reservoir is fully filled to its intended level.

- Effectiveness of Upstream Wall and Cutoff Wall: A wall positioned at the upstream heel has a lesser effect on seepage reduction compared to a cutoff wall below the core. The inclusion of a watertight curtain significantly minimizes seepage through both the dam's body and its foundation, highlighting its importance from the initial stages of dam construction to control leakage effectively.

- Upstream Concrete Cover and Clay Blanket: The performance of the upstream concrete cover and the clay blanket at the reservoir bottom must be evaluated concurrently to maximize their effectiveness in seepage control.

- Combined Seepage Control Measures: The combined application of a clay blanket and an upstream concrete cover demonstrated the highest efficiency in reducing seepage. This synergy should have been considered in the initial design phase to optimize dam performance.

- ❖ Optimized Seepage Control Strategy: The final optimized seepage control model resulted in a leakage rate of less than 1.5 liters per second, which is approximately one-third of the initial model's seepage discharge. This significant reduction confirms the effectiveness of the proposed seepage mitigation measures.

**Declaration of Interest**

We declare that there are no conflicts of interest regarding this manuscript titled "Seepage Analysis and Control of Sahand Rockfill Dam Drainage Using Instrumental Data."

Parvaneh Nikrou

Department of Geography, The University of Alabama, Tuscaloosa, AL 35487, USA

pnikrou@crimson.ua.edu

## 6. Refrences